\documentclass[a4paper,11pt]{article}
\usepackage{hyperref}
\hypersetup{
	colorlinks=true,
	linktoc=all,
	linkcolor=black,     
	citecolor=blue,     
}
\usepackage{tikz}
\usepackage{graphicx}
\usepackage{latexsym, bm}
\usepackage{multicol}
\usepackage{indentfirst}
\usepackage{amssymb,amsfonts}
\usepackage{amsmath}
\usepackage{setspace}
\newtheorem{theorem}{Theorem}[section]
\newtheorem{lemma}[theorem]{Lemma}

\newtheorem{definition}[theorem]{Definition}

\newtheorem{problem}[theorem]{Problem}

\newtheorem{conjecture}[theorem]{Conjecture}

\newcommand{\ignore}[1]{}

\begin{document}
\begin{spacing}{0.95}
\date{}

\title{On a conjecture of Conlon, Fox and Wigderson}
\author{Chunchao Fan,\footnote{Center for Discrete Mathematics, Fuzhou University, Fuzhou, 350108 P.~R.~China. Email: {\tt 1807951575@qq.com}.}\;\;\; \; Qizhong Lin,\footnote{Corresponding author. Center for Discrete Mathematics, Fuzhou University, Fuzhou, 350108 P.~R.~China. Email: {\tt linqizhong@fzu.edu.cn}. Supported in part by National Key R\&D Program of China (Grant No. 2023YFA1010202), NSFC (No.\ 12171088, 12226401) and NSFFJ (No. 2022J02018).}\;\;\; \; Yuanhui Yan \footnote{Center for Discrete Mathematics, Fuzhou University, Fuzhou, 350108 P.~R.~China. Email: {\tt 1366460049@qq.com}.}
}

\maketitle
\begin{abstract}
For graphs $G$ and $H$, the Ramsey number $r(G,H)$ is the smallest positive integer $N$ such that any red/blue edge coloring of the complete graph $K_N$ contains either a red $G$ or a blue $H$.
A book $B_n$ is a graph consisting of $n$ triangles all sharing a common edge.

Recently, Conlon, Fox and Wigderson conjectured that  for any $0<\alpha<1$, the random lower bound $r(B_{\lceil\alpha n\rceil},B_n)\ge (\sqrt{\alpha}+1)^2n+o(n)$ is not tight. In other words, there exists some constant $\beta>(\sqrt{\alpha}+1)^2$ such that $r(B_{\lceil\alpha n\rceil},B_n)\ge \beta n$ for all sufficiently large $n$. This conjecture holds for every $\alpha< 1/6$ by a result of Nikiforov and Rousseau from 2005, which says that in this range $r(B_{\lceil\alpha n\rceil},B_n)=2n+3$ for all sufficiently large $n$.

We disprove the conjecture of Conlon, Fox and Wigderson. Indeed, we show that the random lower bound is asymptotically tight for every $1/4\leq \alpha\leq 1$.
Moreover, we show that for any $1/6\leq \alpha\le 1/4$ and large $n$,
$r(B_{\lceil\alpha n\rceil}, B_n)\le\left(\frac 32+3\alpha\right) n+o(n)$,
where the inequality is asymptotically tight when $\alpha=1/6$ or $1/4$.
We also give a lower bound of $r(B_{\lceil\alpha n\rceil}, B_n)$ for $1/6\le\alpha< \frac{52-16\sqrt{3}}{121}\approx0.2007$, showing that the random lower bound is not tight, i.e., the conjecture of Conlon, Fox and Wigderson holds in this interval.


\medskip

\end{abstract}

\section{Introduction}

For graphs $G$ and $H$, the Ramsey number $r(G,H)$ is the smallest positive integer $N$ such that any red/blue edge coloring of the complete graph $K_N$ contains either a red $G$ or a blue $H$. A central problem in extremal graph theory is to determine the Ramsey number $r(G,H)$.

Let $B_{n}^{(k)}$ be the book graph consisting of $n$ copies of $K_{k+1}$, all sharing a common $K_k$. We always call $n$ the size of the book. When $k=2$, we write $B_{n}$ instead of $B_{n}^{(2)}$. Books have attracted a great deal of attention in graph Ramsey theory, see, for example, the  recent breakthrough of Campos, Griffiths, Morris and Sahasrabudhe \cite{cgms}.

For the diagonal Ramsey number of books, Erd\H{o}s, Faudree, Rousseau and Schelp \cite{efrs} and independently Thomason \cite{tho} proved that $(2^k+o(1))n\le r(B_n^{(k)},B_n^{(k)})\le4^kn.$ In particular, Rousseau and Sheehan \cite{R-S} showed that $r(B_n,B_n)=4n+2$ if $4n+1$ is a prime power.
Recently, a breakthrough result of Conlon \cite{Conlon} established that for each $k\ge2$,
\begin{align*}
r(B_n^{(k)},B_n^{(k)})=2^kn+o(n),
\end{align*}
 which confirms a conjecture of Thomason \cite{tho} asymptotically and also gives an answer to a problem proposed by Erd\H{o}s, Faudree, Rousseau and Schelp \cite{efrs}. The error term $o(n)$ of the upper bound has been improved to $O\big(\frac{n}{(\log\log\log n)^{1/25}}\big)$ by Conlon, Fox and Wigderson \cite{c-f} using a different method,  in particular avoiding the use of the full regularity lemma.

For the off-diagonal Ramsey number of books, a simple lower bound is as follows: for any $k,m,n\in \mathbb{N}$ with $m\le n$,
\begin{align}\label{off-bk-l}
r(B_{m}^{(k)},B_n^{(k)})\ge k(n+k-1)+1.
\end{align}
Indeed, this follows from an observation of Chv\'{a}tal and Harary \cite{c-h} which states that  if $H$ is connected, then $r(G,H)\ge(\chi(G)-1)(|V(H)|-1)+1$. We say that $H$ is $G$-good if this inequality is tight. Burr and Erd\H{o}s \cite{be} initiated the study of such Ramsey goodness problems; the reader is referred to the survey of Conlon, Fox and Sudakov \cite{cfs-15} for a detailed history of the area. Among these results, Nikiforov and Rousseau \cite{nr-to appear} obtained extremely general results. In particular, they proved the following theorem; see Fox, He and Wigderson \cite{fhw} for a new proof avoiding the application of the regularity lemma.

\begin{theorem}[Nikiforov and Rousseau \cite{nr-to appear}]\label{n-r09}
For every $k\ge2$,  there exists some $\alpha_0\in(0,1)$ such that, for any $0<\alpha\le \alpha_0$ and sufficiently large $n$,
\[
r(B_{\lceil\alpha n\rceil}^{(k)},B_n^{(k)})= k(n+k-1)+1.
\]
\end{theorem}

Since we are concerned with the asymptotic behavior of the Ramsey number, we always omit the ceiling and floor signs henceforth.

\medskip
We know that $\alpha_0$ in Theorem \ref{n-r09} is always small. In fact, if $\alpha$ is  sufficiently far from $0$, then the situation of the lower bound is much different. As pointed out by Conlon, Fox and Wigderson \cite{cfw-2}, we can get {\em a random lower bound} for $r(B_{\alpha n}^{(k)},B_n^{(k)})$ as follows (one can also see \cite{C-L-Y2021} for the special case of $k=2$). Indeed, for any $k\in \mathbb{N}$ and $0< \alpha\le 1$, let $p=\frac1{\alpha^{1/k}+1}$ and  $N=(p^{-k}-o(1))n$. Then we randomly and  independently color every edge of $K_N$ blue with probability $p$ and
red with probability $1-p$.  A standard application of the Chernoff bound yields that the probability containing a blue $B_n^{(k)}$ or a red $B_{\alpha n}^{(k)}$ is $o(1)$. Therefore, for any $k\in \mathbb{N}$ and $0< \alpha\le 1$,
\[
r(B_{\alpha n}^{(k)},B_n^{(k)})\ge (\alpha^{1/k}+1)^kn-o(n).
\]
A simple calculation implies that for large $k$, if $\alpha>((1+o(1))\frac{\log k}k)^k$, then the above lower bound is much larger than that of (\ref{off-bk-l}), where the logarithm is to base $e$.

Furthermore, Conlon, Fox and Wigderson \cite{cfw-2}  show that the random lower bound becomes asymptotically tight at this point.

\begin{theorem}[Conlon, Fox and Wigderson \cite{cfw-2}]\label{cfw-2}
For every $k\ge2$,  there exists some $\alpha_1=\alpha_1(k)\in(0,1]$ such that, for any fixed $\alpha_1\le\alpha\le1$,
\[
r(B_{\alpha n}^{(k)},B_n^{(k)})= (\alpha^{1/k}+1)^kn+o(n).
\]
Moreover,  one may take $\alpha_1(k)=((1+o(1))\frac{\log k}k)^k$.
\end{theorem}

From Theorem \ref{n-r09} and Theorem \ref{cfw-2}, we know that for every $k\ge2$, there exist $\alpha_0,\alpha_1\in(0,1]$ such that the following holds:

\medskip
(i) if $0<\alpha\le \alpha_0$, then $r(B_{\alpha n}^{(k)},B_n^{(k)})= k(n+k-1)+1$;

(ii) if $\alpha_1\le\alpha\le1$, then $r(B_{\alpha n}^{(k)},B_n^{(k)})= (\alpha^{1/k}+1)^kn+o(n).$

\medskip
However, the values of $\alpha_0$ and $\alpha_1$ are not very clear in general. 
Moreover, although Theorem \ref{cfw-2}  shows that the random lower bound becomes asymptotically tight when $k\to \infty$, it does not say anything non-trivial in the simplest case $k=2$ except when $\alpha=1$.

Nikiforov and Rousseau \cite{N-R2} proved that
\begin{align}\label{n-r-smb}
\text{for any fixed $\alpha< 1/6$ and all large $n$, $r(B_{\alpha n},B_n)=2n+3$,}
\end{align}
which implies the Ramsey goodness of large books.
Moreover, the constant $1/6$ is asymptotically best possible, i.e., for any $\alpha>1/6$ and all large $n$,
\[r(B_{\alpha n},B_n)>2n+3.\]
Recently, the second author together with Chen and You \cite{C-L-Y2021} proved that for all large $m\le n$, $r(B_m, B_n)\le 2(m+n)+o(n).$ However, the behavior of $r(B_{\alpha n},B_n)$ is still not well-understood for $1/6\le\alpha<1$.

More recently, Conlon, Fox and Wigderson \cite[Conjecture 6.1]{cfw-2} proposed the following conjecture, which would imply that the random lower bound $r(B_{\alpha n},B_n)\ge (\sqrt{\alpha}+1)^2n-o(n)$ is not tight for any $\alpha<1$.
\begin{conjecture}[Conlon, Fox and Wigderson \cite{cfw-2}]\label{cj-cfw}
For every $\alpha<1$, the random lower bound for $r(B_{\alpha n},B_n)$ is not tight. In other words, there exists some $\beta>(\sqrt{\alpha}+1)^2$ such that $r(B_{\alpha n},B_n)\ge \beta n$ for all $n$ sufficiently large.
\end{conjecture}

Note that Conjecture \ref{cj-cfw} holds for any $\alpha\le1/6$, by (\ref{n-r-smb}). However, in this paper, we disprove Conjecture \ref{cj-cfw} by showing that the random lower bound is asymptotically tight for any $1/4\leq \alpha\leq 1$.
\ignore{\medskip
An early and strong problem proposed by Rousseau and Sheehan \cite{R-S} is as follows.
\begin{problem}[Rousseau and Sheehan \cite{R-S}]\label{r-s}
Do there exist constants $a$ and $b$ such that, for all values of $m$ and $n$,
$a\le r(B_m,B_n)-2(m+n)\le b.$
\end{problem}

{\em Remark.} Recently, the second author together with Chen and You \cite{C-L-Y2021} proved that for all large $m\le n$, $r(B_m, B_n)\le 2(m+n)+o(n).$ However, we will see the answer to Problem \ref{r-s} is negative.
}\medskip
\begin{theorem}\label{main2}
For any fixed $1/4\leq \alpha\leq 1$ and sufficiently large $n$,
$$r(B_{\alpha n}, B_n)=(\sqrt{\alpha}+1)^2 n+o(n).$$
\end{theorem}

Moreover, we give an upper bound of $r(B_{\alpha n}, B_n)$ for $1/6\leq\alpha\le 1/4$.
\begin{theorem}\label{upbound}
For any fixed $1/6\leq \alpha\le 1/4$ and sufficiently large $n$,
$$r(B_{\alpha n}, B_n)\le\left(3/2+3\alpha\right) n+o(n).$$
The inequality is asymptotically tight when $\alpha=1/6$ or $1/4$.
\end{theorem}

We also have the following lower bound.

\begin{theorem}\label{lbound}
	Let $1/6\le\alpha\le \frac{52-16\sqrt{3}}{121}\approx0.2007$ be fixed, and let $p=\frac{1-\sqrt{\alpha(3-2\alpha)}}{1-2\alpha}$. Then for all sufficiently large $n$,
	$$r(B_{\alpha n}, B_n)\geq\frac{3n}{1+2p^2}-o(n).$$
The inequality is asymptotically tight when $\alpha=1/6$.
\end{theorem}

{\em Remark 1.}  Since $\frac{3}{1+2p^2}>(\sqrt{\alpha}+1)^2$ for any $1/6\le\alpha< \frac{52-16\sqrt{3}}{121}$, we have that Conjecture \ref{cj-cfw} holds in this interval.

\medskip
\noindent
{\bf  Notation:} For a book $B_n$, we refer to the common edge as the base of the book $B_n$. For a graph $G=(V,E)$ with vertex set $V$ and edge set $E$, we write $bk_G$ for the size of the largest book in a graph $G$. Let $uv$ denote an edge of $G$. For $X \subseteq V$, $e(X)$ is the number of edges in $X$. For two disjoint subsets $X,Y\subseteq V$, $e_G(X,Y)$ denotes the number of edges between $X$ and $Y$. The neighborhood of a vertex $v$ in $U\subseteq V$ is denoted by $N_G(v,U)$, and $\deg_G(v,U)=|N_G(v,U)|$ and the degree of $v$ in $G$ is $\deg_G(v)=|N_G(v,V)|$.
Let $X \sqcup Y$ denote the disjoint union of $X$ and $Y$. Let $[n]=\{1,2,\dots,n\}$, and $[m,n]=\{m,m+1,\dots,n\}$.
We always delete the subscriptions if there is no confusion from the context.

\medskip
The rest of the paper is organized as follows. In Section \ref{pre}, we will collect several useful lemmas. In Sections \ref{pf-1} and \ref{pf-3}, we shall present the proofs of Theorems \ref{main2}, \ref{upbound} and \ref{lbound}.  Finally, we will have some discussion in Section \ref{clu}.

\section{Preliminaries}\label{pre}

The proofs rely on the regularity method \cite{sze78}, which is a powerful tool in extremal graph theory. There are many important applications of the regularity lemma, and we refer the reader to the nice surveys \cite{kss,ko-sim,rs} and many other recent references.

Let $G=(V,E)$ be a graph. For two vertex sets $A,B\subseteq V(G)$, we call $d(A,B)=\frac{e(A,B)}{|A||B|}$ the density of  the pair $(A,B)$.
Let $\varepsilon>0$, a pair $(A,B)$ is said to be \emph{$\varepsilon$-regular} if $|d(A,B)-d(X,Y)|\le \varepsilon$ for every $X\subseteq A$, $Y\subseteq B$ with $|X|\ge \varepsilon |A|$ and $|Y|\ge \varepsilon |B|$. Also,  a subset $A$ is said to be $\varepsilon$-regular if the pair $(A,A)$ is $\varepsilon$-regular.

Given a graph $G$, an equitable $\varepsilon$-regular partition $V(G)=\sqcup_{i=1}^k V_i$ of $G$ is a partition of $V(G)$ such that (i) $||V_i|-|V_j||\le 1$ for all distinct $i$ and $j$; (ii) each $V_{i}$ is $\varepsilon$-regular;  and (iii) for every $1\le i\le k$,  there are at most $\varepsilon k$ values $1\le j\le k$ such that the pair $(V_i, V_j)$ is not $\varepsilon$-regular.

When establishing the asymptotic order of $r(B^{(k)}_n,B^{(k)}_n)$, Conlon \cite[Lemma 3]{Conlon} applied a refined version of the regularity lemma which guarantees a regular subset in each part of the partition for any graph.
We will use the following refined regularity lemma due to Conlon, Fox and Wigderson \cite[Lemma 2.1]{c-f}, which is a strengthening of that due to Conlon \cite{Conlon} and the usual version of Szemer\'{e}di's regularity lemma \cite{sze78}. 
\begin{lemma}[Conlon, Fox and Wigderson \cite{c-f}]\label{regular}
For every $\varepsilon>0$ and $M_0\in \mathbb{N}$, there is some $M=M(\varepsilon, M_0)>M_0$ such that for every graph $G$, there is an equitable $\varepsilon$-regular partition $V(G)=\sqcup_{i=1}^k V_i$ where $M_0\le k\le M$.
\end{lemma}

\ignore{
The following property of regular pairs is immediate from its definition, see e.g. \cite[Fact 1.3]{ko-sim}.

\begin{lemma}\label{pre}
Let $0<\varepsilon<d<1$ and $(A,B)$ be an $\varepsilon$-regular pair of density $d$. If $Y\subseteq B$ with $|Y|\ge \varepsilon |B|$, then there is a subset $A'\subseteq A$ with $|A'|\ge (1-\varepsilon)|A|$ such that each vertex in $A'$ is adjacent to at least $(d-\varepsilon)|Y|$ vertices in $Y$.
\end{lemma}
}

We will use the following version of the counting lemma, proved by Conlon \cite[Lemma 5]{Conlon}. (A similar counting lemma was proved by Nikiforov and Rousseau \cite[Corollary 11]{N-R2}, but they required all of the clusters to be different.) For the general local counting lemma, see R\"{o}dl and Schacht \cite[Theorem 18]{rs}. 

\begin{lemma}[Conlon \cite{Conlon}]\label{count}
For any $\delta>0$, there is $\varepsilon>0$ such that if $U_{1},U_{2}$  (not necessarily distinct), $W_{1},\ldots,W_{l}$ are vertex sets with $(U_{1},U_{2})$ $\varepsilon$-regular of density at least $\delta$ and $(U_{i},W_{j})$ $\varepsilon$-regular of density $d_{ij}$ for all $i\in [2]$ and $j\in [l]$, then there exists an edge $u_1 u_2$ where $u_1\in U_1$ and $u_2\in U_2$ such that $u_1$ and $u_2$ have at least $\sum^{l}_{j=1}(d_{1j}d_{2j}-\delta)|W_{j}|$ common neighbors in $\sqcup_{j=1}^l W_j$.
\end{lemma}


\medskip
We will also use the following well-known Tur\'an's bound.
\begin{lemma} [Tur\'an \cite{turan}]\label{turan}
For any graph $G$ of order $n$ with average degree $d$, the independence number
$\alpha(G)$ is at least $\frac{n}{1+d}.$
\end{lemma}

\section{Proofs of Theorem \ref{main2} and Theorem \ref{upbound}}\label{pf-1}

In the following, for a red/blue edge coloring of $K_N$, we always use $R/B$ to denote the subgraph induced by all red/blue edges.
For every $i,j, (i \neq j$), let $d_{ij}$ be the red density of the pair $(V_i, V_j)$, i.e.,
$d_{ij}=\frac{e_R(V_i, V_j)}{|V_i||V_j|}$.

We will use the following definition introduced by Conlon, Fox and Wigderson \cite{cfw-2}.

\begin{definition}
Fix parameters $\ell\in\mathbb{N}$ and $\varepsilon,\gamma\in (0,1)$ and suppose that we are given a red/blue coloring of $E(K_N)$. Then an $\ell$-tuple of pairwise disjoint vertex sets $V_1,\ldots, V_\ell\subseteq V(K_N)$ is called an $(\ell,\varepsilon,\gamma)$-red-blocked configuration if the following properties are satisfied:

\medskip
1. Each $V_i$ is $\varepsilon$-regular with itself,
\medskip

2. Each $V_i$ has internal red density at least $\gamma$, and
\medskip

3. For all $i\neq j$, the pair $(V_i,V_j)$ is $\varepsilon$-regular and has blue density at least $\gamma$.
\medskip

\noindent
Similarly, we say that $V_1,\ldots, V_\ell$ is an $(\ell,\varepsilon,\gamma)$-blue-blocked configuration if properties (1-3) hold, but the roles of red and blue interchanged.
\end{definition}


We first have a specific structure as in the following lemma, where the proof relies on the refined regularity lemma due to Conlon, Fox and Wigderson \cite{c-f} together with the idea from Nikiforov and Rousseau \cite{N-R2}.

\begin{lemma}\label{structure}
Let $1/6\leq\alpha\leq 1$, $\varepsilon,\gamma\in (0,1)$ where $\varepsilon$ is sufficiently small in terms of $\gamma$, and $n$ is a large integer. Consider a red/blue edge coloring of $K_N$ and an equitable $\varepsilon$-regular  partition of $V(R)$ guaranteed by Lemma \ref{regular} with $N=(x+y\alpha)n+o(n)$, where $x+y=xy$, $1\leq x\leq y\leq 2x$. If $R$ is $B_n$-free and $B$ is $B_{\alpha n}$-free, then the following two properties hold:

\medskip
(1) there exists no $(3,\varepsilon,\gamma)$-red-blocked configuration;

(2) there exists no $(2,\varepsilon,\gamma)$-blue-blocked configuration.
\end{lemma}
{\bf Proof.}
Let $N=(x+y\alpha+\eta)n$ where $\eta>0$ is sufficiently small. Consider a red/blue edge coloring of $K_N$ on vertex set $[N]$. We assume that $\gamma>0$ is taken sufficiently small in terms of $\eta$. Set $\delta=\gamma^2/2$.
Note that $\varepsilon$ is sufficiently small in terms of $\delta$, since $\varepsilon$ is sufficiently small in terms of $\gamma$.

By Lemma \ref{regular}, there is an equitable $\varepsilon$-regular  partition $[N]=\sqcup_{i=1}^k V_i$ for the red graph $R$, i.e., (i) $||V_i|-|V_j||\le 1$ for all distinct $i$ and $j$; (ii) each part $V_{i}$ is $\varepsilon$-regular;  and (iii) for every $1\le i\le k$,  there are at most $\varepsilon k$ values $1\le j\le k$ such that the pair $(V_i, V_j)$ is not $\varepsilon$-regular. Because the colors are complementary, the same conclusion holds for the blue graph.
For convenience, we may assume $|V_i|=N/k=:t$ for all $i\in[k]$.
By the assumption that $R$ is $B_n$-free and $B$ is $B_{\alpha n}$-free, we have
\begin{align}
bk_R&<\frac{1}{x+y\alpha+\eta}N=\frac{1}{x+y\alpha+\eta}kt\le \left(\frac{1}{x+y\alpha}-\gamma\right)kt,\label{eq-3}\\
bk_B&<\frac{\alpha}{x+y\alpha+\eta}N=\frac{\alpha}{x+y\alpha+\eta}kt\le \alpha\left(\frac{1}{x+y\alpha}-\gamma\right)kt.\label{eq-4}
\end{align}

Now we shall prove that there exists no $(3,\varepsilon,\gamma)$-red-blocked configuration.
On the contrary, we may assume that $V_1,V_2,V_3$ is a $(3,\varepsilon,\gamma)$-red-blocked configuration. Let $M$ be the set of all $s\in [k]\setminus[3]$ such that every pair $(V_i,V_s)$ for $i\in [3]$ is $\varepsilon$-regular; clearly, $|M|\ge (k-3)-3\varepsilon k\geq(1-4\varepsilon)k$.

We compute the maximum size of the red books with bases in $E(V_i)$ for $i\in [3]$. Since the red density of the pair $(V_i,V_s)$ is $d_{is}$, we apply Lemma \ref{count} to obtain that
\begin{align}\label{eq-5}
bk_R\ge \sum\limits_{s\in M} (d_{is}^2-\delta)t.
\end{align}
On the other hand, applying Lemma \ref{count}, we obtain that the maximum size $S$ of the blue books with bases in $E(V_1,V_2)$ satisfies
\begin{align}
bk_B\ge S\ge \sum\limits_{s\in M}((1-d_{1s})(1-d_{2s})-\delta )t.\label{eq-6}
\end{align}
Considering in turn $(V_1,V_3)$ and $(V_2,V_3)$, we obtain exactly in the same way,
\begin{align}
bk_B\ge \sum\limits_{s\in M}((1-d_{1s})(1-d_{3s})-\delta )t,\label{eq-7}\\
bk_B\ge \sum\limits_{s\in M}((1-d_{2s})(1-d_{3s})-\delta )t.\label{eq-8}
\end{align}
Let
 $$d_s=\sum_{i=1}^3d_{is},\;\;\text{and} \;\;d_0=\frac{1}{|M|}\sum\limits_{s\in M}d_s.$$
 Adding (\ref{eq-6}), (\ref{eq-7}) and (\ref{eq-8}) each multiplied by $y/3$ together with (\ref{eq-5}) for each $i\in [3]$ multiplied by $x/3$, noting $2x\ge y$, and applying Cauchy's inequality to the double sum we obtain
\begin{align*}
y\cdot bk_B+x\cdot bk_R \geq & \frac{y}{3}\sum\limits_{s\in M}\left(3-2d_s+\sum\limits_{1\le i<j\le3}d_{is}d_{js}-3\delta\right)t+\frac{x}{3}\sum_{s\in M}\left(\sum_{i=1}^3d_{is}^2-3\delta\right)t
\\ \ge & \sum_{s\in M}\left(\frac{y}{3}\left(3-2d_s+\frac{d_s^2}{2}-\frac{1}{2}\sum_{i=1}^3d_{is}^2\right)
+\frac{x}{3}\sum_{i=1}^3d_{is}^2\right)t-(x+y)\delta kt
\\ \geq & \sum_{s\in M}\left(\frac{y}{3}\left(3-2d_s+\frac{d_s^2}{2}\right)+\left(\frac{x}{3}-\frac y6\right)\cdot\frac13{d_s^2}\right)t-(x+y)\delta kt
\\=&\sum_{s\in M}\left(\frac{y}{3}\left(3-2d_s\right)+\frac{x+y}{9}{d_s^2}\right)t-(x+y)\delta kt
\\ \geq& |M|\left(\frac{y}{3}\left(3-2d_0\right)+\frac{x+y}{9}{d_0^2}\right)t-(x+y)\delta kt.
\end{align*}

Therefore, from (\ref{eq-3}) and (\ref{eq-4}), we obtain that
\begin{align*}
|M|\left(\frac{y}{3}\left(3-2d_0\right)+\frac{x+y}{9}{d_0^2}\right)t-(x+y)\delta kt < \left(y\alpha\left(\frac{1}{x+y\alpha}-\gamma\right)+
x\left(\frac{1}{x+y\alpha}-\gamma\right)\right)kt.
\end{align*}
Since $|M|\geq (1-4\varepsilon)k$, we have
$$(1-4\varepsilon)\left(\frac{x+y}{9}d_0^2-\frac{2y}{3}d_0+y\right)-(x+y)\delta< 1-(x+y\alpha)\gamma.$$
Thus, by noting $\delta=\gamma^2/2$, $\varepsilon$ is sufficiently small in terms of $\gamma$, we have
\begin{align*}
\frac{x+y}{9}d_0^2-\frac{2y}{3}d_0+y-1<\frac{1-(x+y\alpha)\gamma+(x+y)\delta}{1-4\varepsilon}-1= \frac{-(x+y\alpha)\gamma+(x+y)\gamma^2/2+4\varepsilon}{1-4\varepsilon}.
\end{align*}
This leads to a contradiction since the right-hand side is negative for sufficiently small $\gamma$ and the left-hand side is non-negative for $\frac{x+y}{9}>0$ and the discriminant of the quadratic form $\Delta=(\frac{2y}{3})^2-\frac{4(x+y)(y-1)}{9}=\frac49(- xy+x+y)=0$ since $x+y=xy$.

It remains to prove that there exists no $(2,\varepsilon,\gamma)$-blue-blocked configuration.
On the contrary, suppose that $V_1,V_2$ is a $(2,\varepsilon,\gamma)$-blue-blocked configuration without loss of generality. Let $M$ be the set of all $s\in [k]\setminus [2]$ such that all pairs $(V_1, V_s)$ and $(V_2, V_s)$ are $\varepsilon$-regular; clearly $|M|\ge (k-2)-2\varepsilon k\geq(1-3\varepsilon)k$.

By a similar argument as aforementioned, we obtain that for $i\in[2]$,
\begin{align}
bk_B\ge \sum_{s\in M}((1-d_{is})^2-\delta)t.\label{eq-9}
\end{align}
Similarly, 
\begin{align}
bk_R\ge \sum\limits_{s\in M}(d_{1s}d_{2s}-\delta)t.\label{eq-10}
\end{align}
Let
 $d_s=\sum_{i=1}^2 d_{is}$, and $d_0=\frac{1}{|M|}\sum_{s\in M}d_s.$ Adding (\ref{eq-10}) multiplied by $x$ together with (\ref{eq-9}) for each $i\in [2]$ multiplied by $y/2$, noting $x\le y$, and applying Cauchy's inequality to the double sum we obtain
\begin{align*}
y\cdot bk_B+x\cdot bk_R \geq& \frac y2\sum_{s\in M}\left(2-2d_s+\sum_{i=1}^2 d_{is}^2-2\delta\right)t+x\sum_{s\in M}\left(d_{1s}d_{2s}-\delta\right)t
\\\ge& \sum_{s\in M}\left(\frac y2(2-2d_s)+\frac {y}{2}\sum_{i=1}^2 d_{is}^2+\frac x2\left(d_s^2-\sum_{i=1}^2 d_{is}^2\right)\right)t-(x+y)\delta kt
\\ \geq& \sum_{s\in M}\left(\frac y2(2-2d_s)+\frac {y+x}{2}\cdot\frac{d_s^2}{2}\right)t-(x+y)\delta kt
\\ \geq& |M|\left(\frac y2(2-2d_0)+\frac {y+x}{4}d_0^2\right)t-(x+y)\delta kt.
\end{align*}

Therefore, from (\ref{eq-3}) and (\ref{eq-4}), we obtain that
\begin{align*}
|M|\left(\frac y2(2-2d_0)+\frac {y+x}{4}d_0^2\right)t-(x+y)\delta kt < \left(y\alpha\left(\frac{1}{x+y\alpha}-\gamma\right)+
x\left(\frac{1}{x+y\alpha}-\gamma\right)\right)kt.
\end{align*}
Since $|M|\geq (1-3\varepsilon)k$, we have
$$(1-3\varepsilon)\left(\frac{y+x}{4}d_0^2-yd_0+y\right)-(x+y)\delta< 1-(x+y\alpha)\gamma.$$
Thus, by noting $\delta=\gamma^2/2$, $\varepsilon$ is sufficiently small in terms of $\gamma$, we have  that
\begin{align*}
\frac{y+x}{4}d_0^2-yd_0+y-1<\frac{1-(x+y\alpha)\gamma+(x+y)\delta}{1-3\varepsilon}-1= \frac{-(x+y\alpha)\gamma+(x+y)\gamma^2/2+3\varepsilon}{1-3\varepsilon}.
\end{align*}
This leads to a contradiction since the right-hand side is negative for sufficiently small $\gamma$ and the left-hand side is non-negative for $\frac{y+x}{4}>0$ and the discriminant of the quadratic form $\Delta=y^2-(y+x)(y-1)=-xy+x+y=0$.
This completes the proof of Lemma \ref{structure}.
\hfill$\Box$

\medskip

Now, we are ready to give proofs for Theorem \ref{main2} and Theorem \ref{upbound}.

\medskip\noindent
{\em Proof sketches of Theorem \ref{main2} and Theorem \ref{upbound}.}
 Consider a red/blue edge coloring of $K_N$ for some suitable $N$. On the contrary, we suppose that $bk_R<n$ and $bk_B<\alpha n$.
Firstly, we apply the refined regularity lemma due to Conlon, Fox and Wigderson \cite{c-f} to obtain an equitable $\varepsilon$-regular partition of $V(R)$ which guarantees the regularity of each cluster with itself. Secondly,  from the assumption that $bk_R<n$ and $bk_B<\alpha n$, Lemma \ref{structure} implies that there are no $(3,\varepsilon,\gamma)$-red-blocked and $(2,\varepsilon,\gamma)$-blue-blocked configurations, which are the bases for subsequent calculations of the corresponding book sizes. According to the densities of clusters, we partition these clusters into red clusters and blue clusters.
For Theorem \ref{main2}, applying the counting lemma due to Conlon \cite{Conlon} to a single red/blue cluster, and combining with Tur\'{a}n's bound on the subgraph $H$ defined on the red clusters, we finally obtain $2bk_R+bk_B\ge 2n + \alpha n$, which leads to a contradiction.
For Theorem \ref{upbound}, based on the lower bounds of the book sizes of $R/B$, we need to consider  the total blue/red densities between red-cluster sets and blue-cluster sets. The situation is more complicated and the computation is much more technical.

The reason we have some improvements is that the specific structure, i.e., there exists no $(3,\varepsilon,\gamma)$-red-blocked and $(2,\varepsilon,\gamma)$-blue-blocked configurations, is more in line with the essence for 2-books than the structure that no $(2,\varepsilon,\gamma)$-red-blocked and $(2,\varepsilon,\gamma)$-blue-blocked configurations used by  Conlon, Fox and Wigderson \cite{cfw-2}.
In particular, the refined regularity lemma by Conlon, Fox and Wigderson \cite{c-f} is a key ingredient of the proofs.

\medskip
\noindent
{\bf Proof of Theorem \ref{main2}.}
Let $1/4\leq\alpha\leq1$ be fixed, and $\eta>0$ is sufficiently small and $n$ is sufficiently large. Let $N=(x+y\alpha+\eta)n$ where $x+y=xy$, $1\leq x\leq y\leq 2x$. Consider a red/blue edge coloring of $K_N$ on vertex set $[N]$. Let $\gamma>0$ be sufficiently small in terms of $\eta$. Set $\delta=\gamma^2/2$. We assume that $\varepsilon>0$ is sufficiently small in terms of $\delta$ and $\gamma$.

By Lemma \ref{regular}, there is an equitable $\varepsilon$-regular  partition $[N]=\sqcup_{i=1}^k V_i$ for the red graph $R$, i.e., (i) $||V_i|-|V_j||\le 1$ for all distinct $i$ and $j$; (ii) each part $V_{i}$ is $\varepsilon$-regular;  and (iii) for every $1\le i\le k$,  there are at most $\varepsilon k$ values $1\le j\le k$ such that the pair $(V_i, V_j)$ is not $\varepsilon$-regular. Because the colors are complementary, the same conclusion holds for the blue graph.
For convenience, we may assume $|V_i|=N/k=:t$ for all $i\in[k]$.
It suffices to show that for all sufficiently large $n$, there exists a red $B_n$ or a blue $B_{\alpha n}$. On the contrary,
\begin{align}
bk_R&<\frac{1}{x+y\alpha+\eta}N=\frac{1}{x+y\alpha+\eta}kt\le \left(\frac{1}{x+y\alpha}-\gamma\right)kt,\label{eq-4-1}\\
bk_B&<\frac{\alpha}{x+y\alpha+\eta}N=\frac{\alpha}{x+y\alpha+\eta}kt\le \alpha\left(\frac{1}{x+y\alpha}-\gamma\right)kt.\label{eq-4-2}
\end{align}

We call a cluster $V_{i}$ \emph{red} if at least half of its internal edges are red and \emph{blue} otherwise.
Clearly, every cluster is either red or blue.

Now we assume that $V_1,\ldots, V_l$ are blue clusters without loss of generality and set $l=\lambda k$, $0\le\lambda\le 1$.
To finish the proof, we will show that (\ref{eq-4-1}) or (\ref{eq-4-2}) is not true.

We first fix a blue cluster $V_i$ and compute the maximum size of the blue books whose bases lie in $E(V_i)$. Let $M_i$ be the set of all $s\in [k]\setminus \{i\}$ such that $(V_i, V_s)$ is an $\varepsilon$-regular pair, and let
$M_{i1} =M_i\cap [l].$
Clearly, $l\ge |M_{i1}|\ge l-1-\varepsilon k\ge(\lambda-2\varepsilon)k$ since $|M_i|\ge (1-\varepsilon)k$.

By Lemma \ref{structure} (2), for every $s\in M_{i1}$, the red density $d_{is}$ of the pair $(V_i, V_s)$ satisfies $d_{is}\le \gamma$, so the blue density of $(V_i, V_s)$ is at least $1-\gamma$. Therefore, by Lemma \ref{count} and noting $\delta=\gamma^2/2$, the maximum size $S$ of the blue books whose bases are in $E(V_i)$ satisfies
\begin{align*}
S\ge \sum_{s\in M_{i1}}\left((1-\gamma)^2-\delta \right)t\ge(1-2\gamma) |M_{i1}|t.
\end{align*}
Then, by noting $|M_{i1}|\ge (\lambda-2\varepsilon)k$, we find that
\begin{align}\label{bk-b-low}
bk_B\ge S \ge  (1-2\gamma)(\lambda-2\varepsilon)kt.
\end{align}

Next we consider the maximum size of the red books whose bases are contained in a red cluster.
Let us define the graph $H$ as follows. The vertex set of $H$ is $[l+1,k]$ and two vertices $i,j\in[l+1,k]$ are joined if and only if the red density $d_{ij}$ of the $\varepsilon$-regular pair $(V_i, V_j)$ satisfies $d_{ij}> 1-\gamma.$

By Lemma \ref{structure} (1), the complement of $H$ is triangle-free, i.e., the independence number of $H$ is at most two, hence Lemma \ref{turan} implies that the average degree of $H$ is at least $(k-l)/2-1$. So there exists a vertex $i\in[l+1,k]$ such that the degree of $i$ in $H$ is at least $(k-l)/2-1$. Let $N_i$ be the set of all $s\in [k]\setminus \{i\}$ such that $(V_i, V_s)$ is an $\varepsilon$-regular pair, and let
$
N_{i1}=N_i\cap N_H(i).$
Clearly, $k-l\ge |N_{i1}|\ge (k-l)/2-1-\varepsilon k$ since $|N_i|\geq (1-\varepsilon )k$.

By Lemma \ref{count} and $\delta=\gamma^2/2$, the maximum size $S$ of the red books whose bases lie in $E(V_i)$ satisfies
\begin{align*}
S\ge \sum_{s\in N_{i1}}\left((1-\gamma)^2-\delta \right)t\geq(1-2\gamma) |N_{i1}|t.
\end{align*}
Then, by noting $|N_{i1}|\ge (k-l)/2-1-\varepsilon k\geq\left(\frac{1-\lambda}{2}-2\varepsilon\right)k$, we find that
\begin{align}\label{bk-r-low}
bk_R\ge S \ge (1-2\gamma)\left(\frac{1-\lambda}{2}-2\varepsilon\right)kt.
\end{align}
Adding (\ref{bk-b-low}) and (\ref{bk-r-low}) multiplied by $2$, we obtain
\begin{align*}
2bk_R+bk_B\geq & (1-2\gamma)\left(1-\lambda-4\varepsilon\right)kt+(1-2\gamma)(\lambda-2\varepsilon)kt
\geq (1-2\gamma-6\varepsilon) kt
\end{align*}

To prove Theorem \ref{main2}, we take $x=\sqrt{\alpha}+1$, $y=1+\frac{1}{\sqrt{\alpha}}$. Since $\varepsilon$ is sufficiently small in terms of $\gamma$ and $\frac{2+\alpha}{x+y\alpha}\leq 1$ for $\frac14\leq\alpha\leq1$, we have $$2bk_R+bk_B\geq\left(\frac{2+\alpha}{x+y\alpha}-2\gamma-\alpha\gamma\right) kt.$$
It follows that either
$2bk_R\geq (\frac{2}{x+y\alpha}-2\gamma) kt$, or
$bk_B\geq (\frac{\alpha}{x+y\alpha}-\alpha\gamma) kt,$
which contradicts (\ref{eq-4-1}) or  (\ref{eq-4-2}), respectively.
The proof of Theorem \ref{main2} is complete.
\hfill$\Box$

\medskip
\noindent
{\bf Proof of Theorem \ref{upbound}.}
Let $1/6\leq\alpha\leq1/4$ be fixed, and $\eta>0$ is sufficiently small and $n$ is sufficiently large. Let $N=(x+y\alpha+\eta)n$, where $x=3/2$ and $y=3$, so $x+y=xy$. Consider a red/blue edge coloring of $K_N$ on vertex set $[N]$. Let $\gamma>0$ be sufficiently small in terms of $\eta$. Set $\delta=\gamma^2/2$, and $\varepsilon>0$ is taken sufficiently small in terms of $\delta$ and $\gamma$.

Similarly, by Lemma \ref{regular}, there is an equitable $\varepsilon$-regular partition $[N]=\sqcup_{i=1}^k V_i$ for the red graph $R$.
For convenience, we may assume $|V_i|=N/k=:t$ for all $i\in[k]$.
It suffices to show that for all sufficiently large $n$, there exists a red $B_n$ or a blue $B_{\alpha n}$. On the contrary,
\begin{align}
bk_R&<\frac{1}{x+y\alpha+\eta}N=\frac{1}{x+y\alpha+\eta}kt\le \left(\frac{1}{x+y\alpha}-\gamma\right)kt,\label{eq-5-1}\\
bk_B&<\frac{\alpha}{x+y\alpha+\eta}N=\frac{\alpha}{x+y\alpha+\eta}kt\le \alpha\left(\frac{1}{x+y\alpha}-20\gamma\right)kt.\label{eq-5-2}
\end{align}

We call a cluster $V_{i}$ \emph{red} if at least half of its internal edges are red and \emph{blue} otherwise.
Clearly, every cluster is either red or blue. Now we assume that $V_1,\ldots, V_l$ are blue clusters without loss of generality and set $l=\lambda k$, $0\le\lambda\le 1$. We first consider a blue cluster $V_i$ and compute the maximum size of the blue books whose bases lie in $E(V_i)$. Let $M_i$ be the set of all $s\in [k]\setminus \{i\}$ such that $(V_i, V_s)$ is an $\varepsilon$-regular pair, and let
\begin{align*}
M_{i1} =M_i\cap [l], \;\;\text{and}\;\;
M_{i2} =M_i\cap [l+1,k].
\end{align*}
Clearly, $l\ge |M_{i1}|\ge l-1-\varepsilon k\geq (\lambda-2\varepsilon) k$ since $|M_i|\ge (1-\varepsilon)k$.

By Lemma \ref{structure} (2), for every $s\in M_{i1}$, the red density $d_{is}$ of the pair $(V_i, V_s)$ satisfies $d_{is}\le \gamma$, so the blue density of $(V_i, V_s)$ is at least $1-\gamma$. Therefore, by Lemma \ref{count} and noting $\delta=\gamma^2/2$, the maximum size $S$ of the blue books whose bases are in $E(V_i)$ satisfies
\begin{align*}
S&\ge \sum_{s\in M_{i1}}\left((1-\gamma)^2-\delta \right)t+\sum_{s\in M_{i2}}((1-d_{is})^2-\delta)t\\&\ge(1-2\gamma) |M_{i1}|t+\sum_{s\in M_{i2}}(1-d_{is})^2 t -\delta|M_{i2}|t.
\end{align*}
Then, by noting $\varepsilon$ is sufficiently small in terms of $\delta$, $|M_{i1}|\ge (\lambda-2\varepsilon)k$ and $|M_{i2}|\leq k$,  and applying Cauchy's inequality, we obtain
\begin{align*}
bk_B\ge S \ge&  (1-2\gamma)(\lambda-2\varepsilon)kt+\sum_{s\in M_{i2}}(1-d_{is})^2 t -\delta kt
\\ \geq& (1-2\gamma)\lambda kt+\frac{1}{|M_{i2}|}\left(\sum_{s\in M_{i2}}(1-d_{is})\right)^2 t -2\delta kt.
\end{align*}
Recall (\ref{eq-5-2}) and $\delta=\gamma^2/2$, we obtain that
\begin{align}\label{range}
\frac{1}{|M_{i2}|}\left(\sum_{s\in M_{i2}}(1-d_{is})\right)^2 t
<\left(\frac{\alpha}{x+y\alpha}-20\gamma-(1-2\gamma)\lambda+2\delta\right)kt
\leq\left(\frac{\alpha}{x+y\alpha}-\lambda-15\gamma\right)kt.
\end{align}
We may assume $\lambda\leq\frac{\alpha}{x+y\alpha}< 1/2$, otherwise the right-hand side of (\ref{range}) is negative, which is not possible since the left-hand side is non-negative.
Since $|M_{i2}|\leq k-l=(1-\lambda)k$, we have
$$
\sum_{s\in M_{i2}}(1-d_{is}) < k\sqrt{\left(\frac{\alpha}{x+y\alpha}-\lambda-15\gamma\right)(1-\lambda)}.
$$
Summing over all $i\in [l]$ and noting that $l=\lambda k<k/2$, we obtain that the total blue densities of all regular pairs $(V_i,V_s)$ where $i\in[l]$ and $s\in[l+1,k]$ satisfies that
\begin{align}\label{bk-b-low-2}
\sum_{i=1}^l \sum_{s\in M_{i2}}(1-d_{is})
< \lambda k^2\sqrt{\left(\frac{\alpha}{x+y\alpha}-\lambda-15\gamma\right)(1-\lambda)}.
\end{align}

Next we consider the maximum size of the red books whose bases are contained in a red cluster.
Let us define the graph $H$ as follows. The vertices of $H$ are the numbers $[l+1,k]$ and two vertices $i,j$ are joined if and only if the red density $d_{ij}$ of the $\varepsilon$-regular pair $(V_i, V_j)$  satisfies
$d_{ij}>1-\gamma.$

By Lemma \ref{structure} (1), the complement of $H$ is triangle-free, hence Lemma \ref{turan} implies that the average degree of $H$ is at least $(k-l)/2-1$. Thus, recall $l\leq k/2$, if $k$ is sufficiently large, then
\begin{align}
e(H)\ge \frac{k-l}{2}\bigg(\frac{k-l}{2}-1\bigg)\geq\left(\frac{1}{4}-\varepsilon\right)(k-l)^2. \label{eH}
\end{align}
 For any $i\in [l+1,k]$, let $N_i$ be the set of all $s\in [k]\setminus \{i\}$ such that $(V_i, V_s)$ is $\varepsilon$-regular, and let
\begin{align*}
N_{i1}=N_i\cap N_H(i), \;\;\text{and}\;\;
N_{i2}=N_i\cap [l].
\end{align*}
Since $|N_i|\ge (1-\varepsilon)k$, we have
$
\deg_H(i)\ge |N_{i1}|\ge \deg_H(i)-\varepsilon k.
$
Therefore, since $\varepsilon$ is sufficiently small in terms of $\delta$, the maximum size $S_i$ of the red books whose bases lie in $E(V_i)$ satisfies
\begin{align*}
S_i\ge& \sum_{s\in N_{i1}}\left((1-\gamma)^2-\delta \right)t+\sum_{s\in N_{i2}}(d_{is}^2-\delta)t
\\ \ge& (1-2\gamma)|N_{i1}|t+\sum_{s\in N_{i2}}d_{is}^2 t-\delta|N_{i2}|t
\\ \ge& (1-2\gamma)\deg_H(i)t+\sum_{s\in N_{i2}}d_{is}^2 t-2\delta kt.
\end{align*}
It follows that
\begin{align*}
\sum_{i=l+1}^k S_i\ge& \sum_{i=l+1}^k (1-2\gamma)\deg_H(i)t+\sum_{i=l+1}^k \sum_{s\in N_{i2}}d_{is}^2 t-2\delta k(k-l)t
\\ \geq& (1-2\gamma)2e(H)t+\sum_{i=l+1}^k\sum_{s\in N_{i2}}d_{is}^2 t-2\delta k(k-l)t.
\end{align*}
Hence by (\ref{eH}) and $\varepsilon$ is sufficiently small in terms of $\delta$, we obtain that
\begin{align*}
\sum_{i=l+1}^k S_i \ge& (1-2\gamma)(1/2-2\varepsilon)(k-l)^2t+\sum_{i=l+1}^k\sum_{s\in N_{i2}}d_{is}^2 t-2\delta k(k-l)t
\\ \ge& \frac{1}{2}(1-2\gamma)(k-l)^2t+\sum_{i=l+1}^k\sum_{s\in N_{i2}}d_{is}^2 t-4\delta k(k-l)t.
\end{align*}
Then by noting that $|N_{i2}|\leq l=\lambda k$ and applying Cauchy's inequality to the double sum, we have
\begin{align*}
bk_R\ge \frac{1}{k-l}\sum_{i=l+1}^k S_i\ge&\frac{1}{2}(1-2\gamma)(k-l)t+\frac{1}{k-l}\sum_{i=l+1}^k\sum_{s\in N_{i2}}d_{is}^2 t-4\delta kt
\\ \geq&\frac{1}{2}(1-2\gamma)(1-\lambda)kt+ \frac{1}{(k-l)l}\sum_{i=l+1}^k\left(\sum_{s\in N_{i2}}d_{is}\right)^2 t-4\delta kt
\\ \geq&\frac{1}{2}(1-2\gamma)(1-\lambda)kt+ \frac{1}{(k-l)^2l}\left(\sum_{i=l+1}^k \sum_{s\in N_{i2}}d_{is}\right)^2 t-4\delta kt.
\end{align*}
So we obtain
\begin{align}\label{bkr}
bk_R\geq\frac{1}{2}(1-2\gamma)(1-\lambda)kt-4\delta kt.
\end{align}
Recall (\ref{eq-5-1}), we obtain that
\begin{align*}
\frac{1}{(k-l)^2l}\left(\sum_{i=l+1}^k \sum_{s\in N_{i2}}d_{is}\right)^2 t
<&\left(\frac{1}{x+y\alpha}-\gamma\right)kt-\frac{1}{2}(1-2\gamma)(1-\lambda)kt +4\delta kt
\\ \leq& \left(\frac{1}{x+y\alpha}-\frac{1-\lambda}{2}+4\delta\right)kt.
\end{align*}
Recall $l=\lambda k$, we obtain that the total red densities  of all regular pairs $(V_i,V_s)$ where $i\in[l+1,k]$ and $s\in[l]$ satisfies that
\begin{align}\label{bk-r-low-2}
\sum_{i=l+1}^k \sum_{s\in N_{i2}} d_{is} <(1-\lambda)k^2\sqrt{\left(\frac{1}{x+y\alpha}-\frac{1-\lambda}{2}+4\delta\right)\lambda}.
\end{align}
Therefore, adding (\ref{bk-b-low-2}) and (\ref{bk-r-low-2}), we obtain
\begin{align*}
\sum_{i=1}^{l}\sum_{s\in M_{i2}}1
&=\sum_{i=1}^l \sum_{s\in M_{i2}}(1-d_{is})+\sum_{i=1}^l \sum_{s\in M_{i2}} d_{is}
=\sum_{i=1}^l \sum_{s\in M_{i2}}(1-d_{is})+\sum_{i=l+1}^k \sum_{s\in N_{i2}} d_{is}
\\&< \lambda k^2\sqrt{\left(\frac{\alpha}{x+y\alpha}-\lambda-15\gamma\right)(1-\lambda)}+
(1-\lambda)k^2\sqrt{\left(\frac{1}{x+y\alpha}-\frac{1-\lambda}{2}+4\delta\right)\lambda}.
\end{align*}
Note that $\sum_{i=1}^{l}\sum_{s\in M_{i2}}1\geq (k-l)l-\varepsilon k^2=((1-\lambda)\lambda-\varepsilon)k^2$, we have
\begin{align*}
((1-\lambda)\lambda-\varepsilon)k^2< \lambda k^2\sqrt{\left(\frac{\alpha}{x+y\alpha}-\lambda-15\gamma\right)(1-\lambda)}
+(1-\lambda)k^2\sqrt{\left(\frac{1}{x+y\alpha}-\frac{1-\lambda}{2}+4\delta\right)\lambda}.
\end{align*}

Suppose $\lambda\leq\eta/10$, then we are done from (\ref{bkr}) that
\[
bk_R\geq\frac{1}{2}(1-2\gamma)(1-\lambda)kt-4\delta kt\ge \left(\frac{1}{2}-\gamma-\frac\lambda2-4\delta\right)\left(\frac32+3\alpha+\eta
\right)n\geq n
\]
by noting $\alpha\ge1/6$, and $\gamma$ and $\delta$ are sufficiently small in terms of $\eta$.
Therefore, we may assume $\lambda>\eta/10$.

Since $\gamma$ is sufficiently small in terms of $\lambda$, $\alpha$, $x$ and $y$, and $\delta=\gamma^2/2$, and $\varepsilon$ is sufficiently small in terms of $\gamma$, we obtain
\begin{align*}
(1-\lambda)\lambda
<& \lambda \sqrt{\left(\frac{\alpha}{x+y\alpha}-\lambda\right)(1-\lambda)}+(1-\lambda)\sqrt{\left(\frac{1}{x+y\alpha}-\frac{1-\lambda}{2}\right)\lambda},
\end{align*}
and consequently,
\begin{align}\label{final}
\sqrt{(1-\lambda)\lambda}< \sqrt{\left(\frac{\alpha}{x+y\alpha}-\lambda\right)\lambda}
+\sqrt{\left(\frac{1}{x+y\alpha}-\frac{1-\lambda}{2}\right)(1-\lambda)}.
\end{align}

Since (\ref{final}) makes sense, we obtain  $1-\frac{2}{x+y\alpha}\leq\lambda\leq\frac{\alpha}{x+y\alpha}$. Recall that $1/6\leq\alpha\leq1/4$, $x=3/2$ and $y=3$, so $\frac{\alpha+{11}/{6}}{x+y\alpha}\leq1$, implying
$\frac{\alpha}{x+y\alpha}
\leq \frac{1}{12}+\frac{11}{12}(1-\frac{2}{x+y\alpha})
\leq \frac{1}{12}+\frac{11}{12}\lambda
= \frac{1}{12}(1-\lambda)+\lambda.$
Thus $\frac{\alpha}{x+y\alpha}-\lambda\leq\frac{1}{12}(1-\lambda)$, and so $\sqrt{(\frac{\alpha}{x+y\alpha}-\lambda)\lambda}\le\frac{1}{2\sqrt{3}}\sqrt{(1-\lambda)\lambda}$.

Moreover, since $\frac{1}{2}-\frac{1}{x+y\alpha}\ge0$, we obtain that
 $$\frac{1}{x+y\alpha}-\frac{1}{2}
\leq\left(2\left(1-\frac{1}{2\sqrt{3}}\right)^2-1\right)\left(\frac{1}{2}-\frac{1}{x+y\alpha}\right)
\leq \left(\left(1-\frac{1}{2\sqrt{3}}\right)^2-\frac12\right)\lambda.$$
Thus $\frac{1}{x+y\alpha}-\frac{1-\lambda}{2}\leq(1-\frac{1}{2\sqrt{3}})^2\lambda$,
and so $\sqrt{(\frac{1}{x+y\alpha}-\frac{1-\lambda}{2})(1-\lambda)}\le(1-\frac{1}{2\sqrt{3}})\sqrt{(1-\lambda)\lambda}$.

Therefore, adding these two terms, the right-hand side of (\ref{final}) is at most $\sqrt{(1-\lambda)\lambda}$, which leads to a contradiction. The proof of Theorem \ref{upbound} is complete.
\hfill$\Box$

\section{Proof of Theorem \ref{lbound}}\label{pf-3}
Let $\frac 16\le\alpha\le \frac{52-16\sqrt{3}}{121}$ be fixed, and $p=\frac{1-\sqrt{\alpha(3-2\alpha)}}{1-2\alpha}$, and $N=(\frac{3}{1+2p^2}-\eta)n$, where $\eta>0$ is sufficiently small and $n$ is sufficiently large. We shall show that for sufficiently large $N$ there exists a partially random red/blue coloring of the edges of $K_{N}$ for which
$$bk_{R}<n,\quad    \text{and}     \quad     bk_{B}<\alpha n.$$

For convenience, assume that $N$ is divisible by 3. Partition $[N]$ into three sets $A_{1},A_{2},A_{3}$, each with $N/3$ vertices, and color the graphs induced by $A_{1},A_{2},A_{3}$ in red. Then edges of the form $uv$ where $u\in A_i,v\in A_{j}$  $(1\le i<j\le3)$ are independently colored red with probability $p$ and blue with probability $1-p$. For $u,v\in A_{i}$, the size of the red book with base $uv$ is a random variable with expected value
\begin{align}\label{srb}
\frac{N}{3}-2+\frac{2N}{3}p^2\leq \frac N3(1+2p^2)= \frac 13\left(\frac{3}{1+2p^2}-\eta\right)(1+2p^2)n=\left(1-\frac{(1+2p^2)\eta}{3}\right)n.
\end{align}

Now suppose $u\in A_{i}$ and  $v\in A_{j}$ where $i\neq j$. If $uv$ is a blue edge, the size of the blue book with base $uv$ is a random variable with expected value
\begin{align*}
\frac{N}{3}(1-p)^2= \frac 13\left(\frac{3}{1+2p^2}-\eta\right)(1-p)^2n=\left(\frac{(1-p)^2}{1+2p^2}-\frac{(1-p)^2\eta}{3}\right)n=
\left(1-\frac{(1+2p^2)\eta}{3}\right)\alpha n,
\end{align*}
by noting that $\alpha=\frac{(1-p)^2}{1+2p^2}$.

By noting (\ref{srb}), if $uv$ is a red edge, then the size of the red book with base $uv$ is a random variable with expected value
\begin{align*}
\frac N3 p^2+\left(\frac{2N}{3}-2\right)p\leq \frac{N}{3}-2+\frac{2N}{3}p^2\leq\left(1-\frac{(1+2p^2)\eta}{3}\right)n.
\end{align*}


We will use the following version of the Chernoff bound \cite[Corollary 2.4 and Theorem 2.8]{jlr}: let $X_1, \ldots, X_t$ be independent random variables taking values in \{0,1\} and let $X=\sum_{i=1}^t X_i$. If $x\geq c\mathbb{E}(X)$ for some $c>1$, then $\mathrm{Pr}(X\geq x)\leq e^{-c'x}$, where $c'=\ln c-1+\frac 1c$.

Plugging in $c=1/(1-\frac{(1+2p^2)\eta}{3})$, since $y=\ln x-1+\frac 1x$ is increasing when $x\geq 1$, we find that $c'>0$. Since $\mathrm{Pr}(X_1\geq n)\leq e^{-c'n}$ and $\mathrm{Pr}(X_2\geq \alpha n)\leq e^{-c'\alpha n}$, where $X_1$ denotes the size of a red book and $X_2$ denotes the size of a blue book, applying a union bound over all edges, we obtain that the probability that there is a red book $B_n$ or a blue book $B_{\alpha n}$ tends to 0 as $n\to\infty$. Thus for large enough $N$ the desired red/blue coloring of the edges of $K_{N}$ exists.\hfill$\Box$

\section{Concluding remarks}\label{clu}
From the result of Nikiforov and Rousseau \cite{N-R2}, we know the exact value of $r(B_{\alpha n}, B_n)$ for $0<\alpha<1/6$; and from Theorem \ref{main2}, we know the asymptotic behavior of $r(B_{\alpha n}, B_n)$ for $1/4\le\alpha\le1$, i.e., the random lower bound $r(B_{\alpha n},B_n)\ge (\sqrt{\alpha}+1)^2n+o(n)$ is asymptotically tight for $1/4\le\alpha\le1$. Moreover, the asymptotic behavior of $r(B_{n}, B_n)$ is already known from a more general result, see \cite{Conlon,c-f}. For the remaining cases, when $1/6\le\alpha\le1/4$, we only know that $r(B_{\alpha n}, B_n)\le (\frac32+3\alpha)n+o(n)$ from Theorem \ref{upbound}. We do not know whether this upper bound is asymptotically tight or not for any $1/6<\alpha<1/4$. Note that  for any $1/6<\alpha<1/4$, $$3/2+3\alpha>(\sqrt{\alpha}+1)^2,$$ therefore, if $r(B_{\alpha n}, B_n)= (\frac32+3\alpha)n+o(n)$ holds in this interval, then it means that Conjecture \ref{cj-cfw} proposed by Conlon, Fox and Wigderson \cite{cfw-2} indeed holds in this interval. In particular, we already know that for any $\frac 16\le\alpha< \frac{52-16\sqrt{3}}{121}\approx0.2007$, Conjecture \ref{cj-cfw} holds from Theorem \ref{lbound}.

\bigskip
{\bf Acknowledgment.}
We are grateful to the anonymous referees for giving invaluable comments and suggestions which improve the presentation of the manuscript greatly.

\end{spacing}

\begin{thebibliography}{99}




\bibitem{be}
S. Burr and P. Erd\H{o}s, Generalizations of a Ramsey-theoretic result of Chv\'{a}tal, {\em J. Graph Theory} 7 (1983), 39--51.

\bibitem{cgms}
M. Campos, S. Griffiths, R. Morris and J. Sahasrabudhe, An exponential improvement for diagonal Ramsey, arXiv:2303.09521v1, 2023.

\bibitem{c-h}
V. Chv\'{a}tal and F. Harary, Generalized Ramsey theory for graphs. III, Small off-diagonal
numbers, {\em Pacific J. Math.} 41 (1972), 335--345.

\bibitem{C-L-Y2021}
X. Chen, Q. Lin and C. You, Ramsey numbers of large books, {\em J. Graph Theory} 101 (2022), 124--133.

\bibitem{Conlon}
D. Conlon, The Ramsey number of books, {\em Adv. Combin.} 3 (2019), 12pp.

\bibitem{cfs-15}
D. Conlon, J. Fox and B. Sudakov, Recent developments in graph Ramsey theory, {\it Surveys in combinatorics 2015}, 49--118. London Math. Soc. Lecture Note Ser.,  424, {\it Cambridge Univ. Press, Cambridge}, 2015.


\bibitem{c-f}
D. Conlon, J. Fox and Y. Wigderson, Ramsey number of books and quasirandomness, {\em Combinatorica} 42 (2022), 309--363.

\bibitem{cfw-2}
D. Conlon, J. Fox and Y. Wigderson, Off-diagonal book Ramsey numbers, {\em Combin. Probab. Comput.} 32 (2023), 516--545.

\bibitem{efrs}
P.~Erd\H{o}s, R.~J.~Faudree, C.~C.~Rousseau and R.~H.~Schelp, The size Ramsey number, {\em Period.~Math.~Hungar.} {9} (1978), 145--161.




\bibitem{fhw}
J. Fox, X. He and Y. Wigderson, Ramsey goodness of books revisited, {\em Adv. Combin.} 4 (2023), 21pp.

\bibitem{kss}
J. Koml\'{o}s, A. Shokoufandeh, M. Simonovits and E. Szemer\'{e}di, The regularity lemma and its
applications in graph theory, {\em Theoretical aspects of computer science (Tehran, 2000)}, Lecture
Notes in Comput. Sci., vol. 2292, Springer, Berlin, 2002, pp. 84--112.

\bibitem{ko-sim}
J. Koml\'{o}s and M. Simonovits, Szemer\'{e}di's regularity lemma and its applications in graph theory.
{\em Combinatorics, Paul Erd\H{o}s is eighty, Vol. 2 (Keszthely, 1993)}, 295--352, Bolyai Soc. Math. Stud., 2, {\em J\'{a}nos Bolyai Math. Soc., Budapest}, 1996.

\bibitem{jlr}
S. Janson, T. {\L}uczak and A. Rucinski, {\em Random graphs}, Wiley-Interscience Series in Discrete Mathematics and Optimization, Wiley-Interscience, New York, 2000.


\bibitem{N-R2} V. Nikiforov and C. C. Rousseau, Book Ramsey numbers I, {\em Random Structures Algorithms} 27 (2005),  379--400.

    \bibitem{nr-to appear}
V. Nikiforov and C. C. Rousseau, Ramsey goodness and beyond, {\em Combinatorica} 29 (2009), 227--262.





\bibitem{rs}
V. R\"{o}dl  and M. Schacht: Regularity lemmas for graphs, in: {\em Fete of Combinatorics
and Computer Science}, Bolyai Soc. Math. Stud. 20, 2010, 287--325.


\bibitem{R-S} C. C. Rousseau and J. Sheehan, On Ramsey numbers for books, {\em J. Graph
Theory} 2 (1978), 77--87.

\bibitem{sze78}
E.~Szemer\'edi, {Regular partitions of graphs}, in: {\em Probl\`emes Combinatories et th\'eorie des graphs, (Colloq. Internat. CNRS, Univ. Orsay, Orsay, 1976)}, pp. 399--401, Colloq. Internat. CNRS, 260, CNRS, Paris, 1978.

\bibitem{tho}
A.~Thomason, On finite Ramsey numbers, {\em European J.~Combin.} {3} (1982), 263--273.
\bibitem{turan}
P. Tur\'an, Eine Extremalaufgabe aus der Graphentheorie [in Hungarian], {\em Math Fiz. Lapok} 48 (1941), 436--452.
\end{thebibliography}
\end{document}